\documentclass[12pt]{article}
\usepackage[cp1251]{inputenc}
\usepackage[russian]{babel}
\usepackage{amssymb,amsthm,amsmath}

\theoremstyle{plain}

\newtheorem{Theorem}{Теорема}
\newtheorem{Statement}{Утверждение}

\newcommand{\intl}{\mathop{\int}\limits}
\newcommand{\liml}{\mathop{\lim}\limits}
\newcommand{\suml}{\mathop{\sum}\limits}

\renewcommand{\le}{\leqslant}
\renewcommand{\ge}{\geqslant}

\def\la{\lambda}
\def\eps{\varepsilon}

\def\N{{\mathbb N}}

\begin{document}

{\bf УДК 517.984}

{\large {\bf \centerline{О скорости равносходимости разложений в
ряды}} {\bf \centerline{по тригонометрической системе и по
собственным функциям}} {\bf \centerline{оператора Штурма--Лиувилля
с потенциалом --- распределением }} }

\bigskip

{\centerline{Садовничая~И.~В.}}

\medskip
В настоящей работе изучается оператор Штурма--Лиувилля
\begin{equation}\label{main}
Ly=l(y)=-\dfrac{d^2y}{dx^2}+q(x)y,
\end{equation}
в пространстве \(L_2[0,\pi]\) с граничными условиями Дирихле
$y(0)=y(\pi)=0$. Предполагается, что потенциал  \(q(x)=u'(x),\)
\(u\in L_2[0,\pi]\). Производная здесь понимается в смысле
распределений. Операторы такого вида были определены в работе
\cite{SSh1}. В работах \cite{SSh1} -- \cite{SSh2} было доказано,
что оператор \(L\) фредгольмов с индексами (0,0) (в случае
вещественного потенциала -- самосопряжен), полуограничен, имеет
чисто дискретный спектр. Некоторые результаты об операторе \(L\),
полученные в \cite{SSh1} -- \cite{SSh3} и необходимые нам в данной
работе, будут приведены ниже.

В статье рассматривается вопрос о равномерной на всем отрезке
\([0,\pi]\) равносходимости разложения функции \(f\) в ряд по
системе собственных и присоединенных функций опрератора \(L\) с ее
разложением в ряд Фурье по системе синусов. Эта задача хорошо
известна в классической теории операторов Штурма--Лиувилля (в
случае, когда потенциал локально суммируем). В монографии
В.~А.~Марченко \cite[\S 3, гл. 1]{Mar} была доказана равномерная
равносходимость в случае, если $f\in L_2[0,\pi]$, а $q$ ---
комплекснозначная суммируемая функция. В 1991 году В.~А.~Ильин
получил результаты в случае, когда $f\in L_1[0,\pi]$,  $q$ ---
комплекснозначная суммируемая функция (см. \cite{Il}).
В.~А.~Винокуров и В.~А.~Садовничий в \cite{VS} доказали теорему о
равносходимости для случая операторов с потенциалом ---
производной функции ограниченной вариации, при этом \(f\in
L_1[0,\pi]\).

Здесь мы покажем, что, если первообразная \(u\) от потенциала ---
комплекснозначная функция из пространства \(L_2[0,\pi]\), то для
любой функции \(f\) из пространства \(L_2[0,\pi]\) имеет место
равномерная на всем отрезке \([0,\pi]\) равносходимость разложений
в ряды по системе синусов и по системе собственных и
присоединенных функций оператора \(L\). Некоторые результаты о
такой равносходимость были анонсированы автором в работе
\cite{Sad}, но подробное доказательство приведено не было.

\vskip 2cm

{\bf \S 1. Предварительные результаты.}

\medskip

 Нам понадобятся некоторые предварительные сведения об
операторе \eqref{main}. Обозначим через $\omega(x,\la)$ решение
дифференциального уравнения $-\omega''+q\omega=\la\omega$ с
начальными условиями $\omega(0,\la)=0$, $\omega^{[1]}(0,\la)=1$
(здесь $\omega^{[1]}=\omega'-u\omega$ --- первая
квазипроизводная). Ясно, что нули целой функции $\omega(\pi,\la)$
совпадают с собственными значениями оператора $L$. В силу теоремы
существования и единственности, геометрическая кратность каждого
собственного значения равна $1$. В случае вещественного потенциала
присоединенные функции отсутствуют, т.е. все собственные значения
являются простыми, но в общем случае это не так. Обозначим через
$\{\mu_k\}_1^\infty$ все нули функции $\omega(\pi,\la)$ \emph{без
учета кратности} (т.е. $\mu_j=\mu_k$ только при $j=k$), причем
нумерацию будем вести в порядке возрастания модуля, а в случае
совпадения модулей --- по возрастанию аргумента, значения которого
выбираются из полуинтервала $(-\pi,\pi]$. Цепочкой из собственной
и присоединенных функций, отвечающей собственному значению
$\mu_k$, называют систему функций
$\{y_k^0,\,y_k^1,\,\dots,\,y_k^{p_k-1}\}$, где
$(L-\mu_kI)y_k^0=0$, $(L-\mu_kI)y_k^j=y_k^{j-1}$ для всех
$j=1,\,\dots,\,p_k-1$. Присоединенные функции определены
неоднозначно и, следуя работе \cite{Ke}, \emph{канонической
цепочкой} назовем любую из цепочек, имеющих максимальную длину.
Длину канонической цепочки (число $p_k$) назовем
\emph{алгебраической кратностью} собственного значения $\mu_k$. В
нашем случае каноническую цепочку можно предъявить явно, что мы
сейчас и проделаем.

Пусть $\mu_k$ --- ноль функции $\omega(\pi,\la)$ кратности $q_k$.
Заметим, что функции $\omega^{(j)}_{\la}(x,\la)$,
$j=1,2,\dots,q_k-1$ удовлетворяют дифференциальным уравнениям
$l(\omega^{(j)}_{\la})=\la\omega^{(j)}_{\la}+\omega^{(j-1)}_\la$,
причем $\omega^{(j)}_{\la}(0,\la)\equiv0$. Кроме того,
$\omega^{(j)}_\la(\pi,\mu_k)=0$, поскольку $\mu_k$ есть ноль
кратности $q_k$. Таким образом, функции
$\omega^{(j)}_\la(x,\mu_k)$, $j=0,1,\dots,q_k-1$ образуют цепочку
из собственной и присоединенных функций, отвечающую собственному
значению $\mu_k$. Как и любая цепочка из собственной и
присоединенных функций, эта система линейно независима, а значит,
порождает подпространство размерности $q_k$. Теперь остается
отметить (см. \cite[Гл. I, п.3]{Na}), что эта цепочка является
канонической, т.е. имеет максимальную длину ($p_k=q_k$).

Нам будет удобно ввести и другую нумерацию собственных и
присоединенных функций.  Под системой собственных и присоединенных
функций $\{y_n(x)\}_{n=1}^\infty$ мы будем везде далее понимать
систему, полученную нормировкой $\|y_n\|_{L_2}=1$ системы
$\bigcup\limits_{k=1}^{\infty}\{\omega^{(j)}_\la(x,\mu_k)\}_{j=0}^{q_k-1}$,
а через $\{\la_n\}_{n=1}^\infty$ обозначим нули функции
$\omega(\pi,\la)$ в порядке возрастания модуля \emph{с учетом
кратности} (тогда $Ly_n=\la_ny_n$ для любого $n\in\N$).

Построим теперь биортогональную систему $\{w_n\}_{n=1}^\infty$ к
системе $\{y_n\}_{n=1}^\infty$. Прежде всего заметим, что при
любом $k$ система $\{z_k^j\}_{j=0}^{q_k-1}$, где
$z_k^j(x)=\overline{\omega^{(j-1)}_\la(x,\mu_k)}$, является
системой из собственной и присоединенных функций для оператора
$L^*=-\frac{d^2}{dx^2}+\overline{q}$, отвечающей собственному
значению $\overline{\mu_k}$. Покажем теперь, как с помощью
конечных линейных комбинаций системы
$\bigcup\limits_{k=1}^{\infty}\{z_k^j\}_{j=0}^{q_k-1}$ построить
биортогональную к  $\{y_n\}_1^\infty$ систему. Заметим, что если
функция $z$ лежит в корневом подпространстве оператора $L^*$,
отвечающем собственному значению $\overline{\mu_l}$, а функция $y$
--- в корневом подпространстве оператора $L$, отвечающем
собственному значению $\mu_k$, где $k\ne l$, то функции $y$ и $z$
ортогональны. Действительно, для собственных функций
$\mu_k(y,z)=(Ly,z)=(y,L^*z)=\mu_l(y,z)$, т.е. $(y,z)=0$. Если
теперь $y$ --- первая присоединенная функция, а $z$ ---
собственная функция, то, учитывая, что ортогональность собственных
функций уже доказана, $\mu_k(y,z)=(Ly,z)=(y,L^*z)=\mu_l(y,z)$.
Дальнейшее очевидно. Таким образом, достаточно построить
биортогональную систему в каждом корневом подпространстве
по--отдельности. В случае простого собственного значения это
легко: $w_n(x)=\overline{y_n(x)}/(y_n(x),\overline{y_n(x)})$, а
для собственных значений кратности $p_k>1$ мы сошлемся на \cite[\S
2]{Sh} или \cite[Гл. 1, \S 3]{Mar}:
$$
w_k^j(x)=-b_k^{p_k-j-1}z_k^{p_k-j-1}+(\overline{\mu_k})^{-1}b_k^{p_k-j-2}z_k^{p_k-j-2}-
2(\overline{\mu_k})^{-2}b_k^{p_k-j-3}z_k^{p_k-j-3}-\dots
-j(-\overline{\mu_k})^{-j}b_k^0z_k^0,
$$
где $b_k^j$ --- некоторые ненулевые числа такие, что разложение
резольвенты $(L-\la I)^{-1}$ в ряд Лорана в окрестности точки
$\mu_k$ имеет вид
$$
(L-\la I)^{-1}=\suml_{s=0}^{p_k-1}\frac{(\cdot,b_k^0z_k^0)y_k^0+
(\cdot,b_k^1z_k^1)y_k^1+\dots+(\cdot,b_k^sz_k^s)y_k^s}{(\la-\mu_k)^{p_k-s}}.
$$
\begin{Statement}\label{st:5}(см. \cite{SSh2}, теорема 2.7.)
Пусть $u\in L_2[0,\pi]$. Тогда система $\{y_n(x)\}_{n=1}^\infty$
собственных и присоединенных функций оператора \(L\) образует
базис Рисса в пространстве $L_2[0,\pi]$.
\end{Statement}
\begin{Statement}\label{st:7}(см. \cite{SSh2}, теорема 3.13.)
Пусть $u\in L_2[0,\pi]$. Тогда, начиная с некоторого номера
$N=N_{u}$, все собственные значения оператора $L$ просты, а для
функций $y_n$ и $w_n$ справедливы асимптотические равенства
\begin{equation}\label{efas}
\begin{array}{c}
y_n(x)=\sqrt{\frac2\pi}\sin nx+\varphi_n(x),\quad
w_n(x)=\sqrt{\frac2\pi}\sin nx+\psi_n(x),\\
y'_n(x)=n\left(\sqrt{\frac{2}{\pi}}\cos nx+\eta_n(x)\right)+
u(x)\left(\sqrt{\frac{2}{\pi}}\sin nx+\varphi_n(x)\right)\quad
n=N,N+1,\dots\ ,
\end{array}
\end{equation}
причем последовательность
$\{\gamma_n\}_{n=N}^\infty=\{\|\varphi_n(x)\|_{C}+
\|\psi_n(x)\|_{C}+\|\eta_n(x)\|_{C}\}_{n=N}^\infty\in l_2$ и ее
норма в этом пространстве ограничена постоянной, зависящей только
от $u$. Кроме этого, при $n\ge N$
$$
\psi_n(x)=\psi_{n,0}(x)+\psi_{n,1}(x)+\psi_{n,2}(x).
$$ Здесь
\begin{equation} \label{efas0}
\psi_{n,0}(x)=\alpha_n \sin nx+\beta_n\cos nx-\intl_0^xu(t)\sin
n(x-2t)dt,
\end{equation}
 где
норма
 последовательности
$\{|\alpha_n|+|\beta_n|\}_{n=N}^\infty$ в пространстве $l_2$
ограничена величиной $C_{u}$;
$$
\psi_{n,1}(x)=\sin
nx\left(-\frac{1}{2n}\intl_0^xu^2(t)\sin(2nt)dt\right)+\cos
nx\left(-\frac{x}{\pi}\intl_0^\pi u(t)\sin(2nt)dt+\right.
$$\begin{equation} \label{efas1} \left.+\frac{x}{2\pi n}\intl_0^\pi
u^2(t)(\cos(2nt)-1)dt-\frac{2x}{\pi}\intl_0^{\pi}\intl_0^t
u(t)u(s)\cos(2nt)\sin(2ns)dsdt+\right.
\end{equation}
$$
\left.+\frac{1}{2n}\intl_0^x
u^2(t)(1-\cos(2nt))dt+2\intl_0^x\intl_0^t u(t)u(s)\cos(2nt)\sin(2ns)dsdt\right),
$$
а последовательность $\{\|\psi_{n,2}(x)\|_C\}_{n=N}^\infty$
принадлежит пространству $l_1$ и ее норма в этом пространстве не
превосходит $C_u$.
\end{Statement}

\vskip 2cm

{\bf \S 2. Основная теорема.}
\medskip

\begin{Theorem}
Рассмотрим оператор \eqref{main}, действующий в пространстве
\(L_2[0,\pi]\), с граничными условиями Дирихле, потенциал которого
удовлетворяет следующим условиям:\,
 \(q(x)=u'(x)\), где комплекснозначная функция \(u\in L_2[0,\pi]\). Пусть \(\{
y_n(x)\}_{n=1}^\infty\) -- нормированная система собственных и
присоединенных функций оператора \(L\), \(\{
w_n(x)\}_{n=1}^\infty\) -- биортогональная к ней система. Для
произвольной функции \(f\in L_2 [0,\pi]\) обозначим
\(c_n=(f(x),w_n(x))\), \(c_{n,0}=\sqrt{2/\pi}(f(x),\sin{nx})\).
Тогда имеет место равномерная на всем отрезке \([0,\pi]\)
равносходимость разложения функции \(f\) в ряд по системе \(\{
y_n(x)\}_{n=1}^\infty\) и по системе синусов:
\begin{equation}
\label{1} \left\|\sum_{n=1}^{m}{c_n
y_n(x)}-\sum_{n=1}^{m}\sqrt{\frac{2}{\pi}}{c_{n,0}\sin
nx}\right\|_{C}\le
C_u\left(\suml_{n\ge\sqrt{m}}|c_{n,0}|^2\right)^{1/2}+\|f\|_{L_2}\upsilon_u(m),
\end{equation}
где $\upsilon_u(m)\to0$ при $m\to +\infty$.
\end{Theorem}
{\bf Замечание.}\quad В работе \cite{Sad1} было доказано, что в
случае более гладкого потенциала, когда $u\in W_2^\theta[0,\pi]$
при $\theta\in(0,1/2)$, выполнена оценка $\upsilon_u(m)\le
C(\|u\|,\theta,\eps)m^{-\theta/2+\eps}$ для любого $\eps>0$.

{\bf Доказательство теоремы 1.}\quad Рассмотрим операторы
\(B_{m,N}: L_2[0,\pi]\to C[0,\pi]\), действующие по правилу
$$
B_{m,N}f(x):=\sum_{n=N}^{m}{c_ny_n(x)}-\sum_{n=N}^{m}\sqrt{\frac{2}{\pi}}{c_{n,0}\sin
nx},
$$ где \(f\in L_2[0,\pi].\) При $N=1$ будем опускать второй
индекс: $B_{m,1}=:B_m$.

Мы зафиксируем индекс $N$, положив $N=N_u$, где $N_u$ ---
натуральное число, существование которого постулируется в
утверждении 2. Отметим, что в случае вещественного потенциала
$N_u=1$.

Очевидно, что для любой функции $f\in L_2[0,\pi]$ справедливо
асимптотическое равенство (мы будем называть его основным
равенством)
$$
B_{m,N}f(x)=\sum_{n=N}^{m}\sqrt{\frac{2}{\pi}}(f(t),\psi_{n,0}(t))\sin
nx+\sum_{n=N}^{m}\sqrt{\frac{2}{\pi}}(f(t),\psi_{n,1}(t))\sin nx+
$$
\begin{equation}
\label{3}
+\sum_{n=N}^{m}\sqrt{\frac{2}{\pi}}(f(t),\psi_{n,2}(t))\sin nx+
\sum_{n=N}^{m}\sqrt{\frac{2}{\pi}}(f(t),\sin
nt)\varphi_n(x)+\sum_{n=N}^{m}(f(t),\psi_n(t))\varphi_n(x).
\end{equation}

Теперь мы хотим показать, что для любой функции $u\in L_2[0,\pi]$
последовательность операторов $\{B_{m,N}\}_{m=N}^\infty$
равномерно ограничена:
\begin{equation} \label{B}
\|B_{m,N}\|_{L_2\to C}\le C_u.
\end{equation}

Оценим каждое из слагаемых в правой части основного равенства
\eqref{3} по отдельности (наиболее трудной здесь будет оценка
первого слагаемого).

\medskip
{\bf  Шаг 1 (оценка первого слагаемого).}
\medskip

{\it Для любой функции $f\in L_2[0,\pi]$ имеет место неравенство
\begin{equation}
\label{0} \left \|\sqrt{\frac{2}{\pi}}\sum_{n=N}^{m}(f(t),
\psi_{n,0}(t))\sin nx\right\|_{C}\le C_u\|f\|_{L_2}.
\end{equation} }

В силу асимптотических формул \eqref{efas0}

\begin{equation*}
\label{4} (f(t),\psi_{n,0}(t))=\overline{\alpha_n}(f(t),\sin
nt)+\overline{\beta_n}(f(t),\cos nt)-
\intl_0^{\pi}f(t)\intl_0^{t}\overline{u(s)}\sin n(t-2s)dsdt.
\end{equation*}
Таким образом, нам необходимо оценить равномерно по всем
$x\in[0,\pi]$ выражение

\begin{equation}
\label{4'} \sum_{n=N}^{m}\left(\overline{\alpha_n}(f(t),\sin
nt)+\overline{\beta_n}(f(t),\cos nt)
-\intl_0^{\pi}f(t)\intl_0^{t}\overline{u(s)}\sin
n(t-2s)dsdt\right)\sin nx
 \end{equation}
 Так как $\|\{\alpha_n\}\|_{l_2}\le C_u$, то
 \begin{equation*}
% \label{4''}
\left
\|\sum_{n=N}^{m}\sqrt{\frac{2}{\pi}}\overline{\alpha_n}(f(t), \sin
nt)\sin nx\right\|_{C}\le C_u\|f\|_{L_2}.
 \end{equation*}
 Аналогичные рассуждения справедливы и для второго слагаемого в
 \eqref{4'}. Наиболее сложными для оценки является третье
 слагаемое.
 Применив формулу разности косинусов, получим
 $$
 \suml_{n=N}^{m}\intl_0^\pi
 f(t)\intl_0^t\overline{u(s)}\sin n(t-2s)dsdt\sin nx=
 \frac 12\suml_{n=N}^{m}\intl_0^\pi
 f(t)\intl_0^t\overline{u(s)}\cos n(t-x-2s)dsdt-
 $$
 $$-\frac 12\suml_{n=N}^{m}\intl_0^\pi
 f(t)\intl_0^t\overline{u(s)}\cos n(t+x-2s)dsdt.
 $$
 Оценим первое слагаемое в правой части последнего соотношения
 (второе слагаемое рассматривается аналогично). Заметим, что
$$
\suml_{n=N}^{m}\intl_0^\pi
 f(t)\intl_0^t\overline{u(s)}\cos n(t-x-2s)dsdt=\intl_0^\pi
 f(t)\intl_0^t\overline{u(s)}D_m(t-x-2s)dsdt-
 $$
 \begin{equation}
 \label{5}-\intl_0^\pi
 f(t)\intl_0^t\overline{u(s)}D_{N-1}(t-x-2s)dsdt,
 \end{equation}
где $D_m(\xi)=1/2+\suml_{n=1}^{m}\cos n\xi$ -- ядро Дирихле. Так
как второе слагаемое в правой части \eqref{5} не зависит от \(m\),
то для него сразу получаем необходимую оценку:
$$\left\|\intl_0^\pi
 f(t)\intl_0^t\overline{u(s)}D_{N-1}(t-x-2s)dsdt\right\|_C\le
 C_u\|f\|_{L_2}.
 $$
 Займемся оценкой первого слагаемого в \eqref{5}. Определим
 оператор \(A_{m,-x}\), действующий в пространстве \(L_2[0,\pi]\)
  по правилу:
 $$
 A_{m,-x}u(t)=\intl_0^t u(s)D_m(t-x-2s)ds
 $$
 и положим по определению $A_m:=A_{m,0}$. Нам необходимо оценить
 нормы операторов $A_{m,-x}$ равномерно по $m\in\N$ и $x\in [0,\pi]$.

 Покажем, что оператор $A_{m,-x}$ унитарно эквивалентен сумме оператора $A_m$ и
 некоторого оператора $\tilde{A}_{m,x}$, норму
 которого в \(L_2[0,\pi]\) оценить легко. Пусть $T_xg(t)=g(t+x)$
 --- оператор сдвига в пространстве \(L_2[0,\pi]\) (считаем, что все
 функции продолжены за отрезок $[0,\pi]$ периодически). Тогда можем
 записать
 $$
 T_{-x}A_{m,-x}T_xu(t)=T_{-x}A_{m,-x}u(x+t)=T_{-x}\intl_0^tu(s+x)D_m(t-x-2s)ds=
 $$
 $$=\intl_0^{t-x}u(s+x)D_m(t-2x-2s)ds=\intl_x^tu(s)D_m(t-2s)ds
 =A_mu(t)-\tilde{A}_{m,x}u(t),
 $$
 где $\tilde{A}_{m,x}u(t):=\intl_0^xu(s)D_m(t-2s)ds$. Несложно
 видеть, что оператор $\tilde{A}_{m,x}$ представляет собой
 композицию оператора срезки $H_xu(t)=\chi_{[0,x]}u(t)$ и
 оператора $S_m$, действующего по правилу
 \begin{equation}
 \label{7}
 S_mv(t)=\suml_{n=1}^m\left(\intl_0^\pi v(s)\cos 2nsds\cos
 nt+\intl_0^\pi v(s)\sin 2nsds\sin nt\right)+\frac 12\intl_0^\pi
 v(s)ds.
 \end{equation}
 Очевидно, что $\|H_x\|_{L_2}\le 1$,
 $\|S_mv(t)\|_{L_2}=\left\|\suml_{n=1}^m(a_n\cos nt+b_n\sin nt)+\frac
 12a_0\right\|_{L_2}^2\le 2\|v\|_{L_2}^2$, где
 $a_n=\intl_0^\pi v(s)\cos 2nsds$,  $b_n=\intl_0^\pi v(s)\sin
 2nsds$ -- коэффициенты Фурье функции $v(t)$. Мы
 показали, что норма оператора $\tilde{A}_{m,x}$ в \(L_2[0,\pi]\)
 не превосходит числа $\sqrt{2}$ при любых значениях $m\in\N$ и $x\in
 [0,\pi]$.

 Теперь нам осталось доказать равномерную по $m\in\N$
 ограниченность оператора $A_m$ в пространстве \(L_2[0,\pi]\).
 Заметим, что
 $$A_mu(t)=\intl_0^tu(s)D_m(t-2s)ds=\intl_0^\pi u(s)D_m(t-2s)ds
 -\intl_t^\pi u(s)D_m(t-2s)ds=S_mu(t)-E_m^*u(t),$$
где оператор $S_m$ определен в \eqref{7}, а $E_m^*$ является
сопряженным к оператору $E_m$, определенному равенством
$E_mu(t)=\intl_0^t u(s)D_m(2t-s)ds$. Таким образом, задача свелась
к проверке равномерной ограниченности оператора $E_m$.

Поскольку $E_m$ --- оператор типа Харди, то его ограниченность
может быть доказана стандартным приемом (см., например, \cite[Гл.
IX, п. 3]{Ha}). Рассмотрим билинейную форму $(E_mf,g)$, где
$f,g\in L_2[0,\pi]$ и докажем ее ограниченность. В силу известной
оценки $|D_m(x)|\le\frac{C}{|x|}$, где $C$ -- некоторая абсолютная
постоянная (см., например, \cite[Гл. 1 \S 32]{Ba}), имеем:
$$|(E_mf,g)|=\left|\intl_0^\pi\intl_0^t
f(s)\overline{g(t)}D_m(2t-s)dsdt\right|\le\intl_0^\pi\intl_0^t |f(s)||g(t)||D_m(2t-s)|dsdt\le$$
$$
\le C\intl_0^\pi\intl_0^t|f(s)||g(t)|\frac{1}{2t-s}dsdt=
C\intl_0^\pi\intl_0^t|f(s)|\frac{1}{\sqrt{2t-s}}\sqrt[4]{\frac
st}|g(t)|\frac{1}{\sqrt{2t-s}}\sqrt[4]{\frac ts}dsdt\le
$$
$$
\le\left(\intl_0^\pi\intl_0^t|f(s)|^2\frac{1}{2t-s}\sqrt{\frac
st}dsdt\right)^{1/2}\left(\intl_0^\pi\intl_0^t|g(t)|^2\frac{1}{2t-s}\sqrt{\frac
ts}dsdt\right)^{1/2}.
$$
Обозначим двойные интегралы в последней части цепочки неравенств
через $I_1$ и $I_2$ соответственно и оценим каждый из них в
отдельности. В интеграле $I_1$ поменяем пределы интегрирования.
Получим $$
I_1=\intl_0^\pi\intl_s^\pi|f(s)|^2\frac{1}{2t-s}\sqrt{\frac
st}dtds=\intl_0^\pi|f(s)|^2\sqrt{s}\intl_s^\pi\frac{d(2\sqrt{t})}{2t-s}ds=
$$
$$=\frac{1}{\sqrt2}
\intl_0^\pi|f(s)|^2\left(\ln\left(\frac{\pi-\sqrt{s/2}}{\pi+\sqrt{s/2}}\right)-2\ln(\sqrt{2}-1)\right)ds,
$$
то есть $\sqrt{I_1}\le (\ln 4)\|f\|_{L_2}$.

Аналогично
$$
I_2=\intl_0^\pi\intl_t^\pi|g(t)|^2\frac{1}{2t-s}\sqrt{\frac
ts}dsdt=\intl_0^\pi|g(t)|^2\sqrt{t}\intl_0^t\frac{d(2\sqrt{s})}{2t-s}ds=
\frac{1}{\sqrt2} \intl_0^\pi|g(t)|^2(2\ln(\sqrt{2}+1))dt,
$$
значит, $\sqrt{I_2}\le C\|g\|_{L_2}$,
$C=\sqrt{\sqrt{2}\ln(\sqrt{2}+1)}$.

Неравенство \eqref{0} полностью доказано.

\medskip
{\bf  Шаг 2 (оценка второго слагаемого в основном равенстве).}
\medskip

{\it Пусть $u\in L_2[0,\pi]$. Тогда для любой функции $f\in
L_2[0,\pi]$ имеет место неравенство
\begin{equation}
\label{0'} \left \|\sqrt{\frac{2}{\pi}}\sum_{n=N}^{m}(f(t),
\psi_{n,1}(t))\sin nx\right\|_{C}\le C_u\|f\|_{L_2}.
\end{equation} }

Заметим, что в силу асимптотических формул \eqref{efas1} нам нужно
рассмотреть выражение

$$
\Psi_n=-\frac{1}{2n}\intl_0^\pi
f(t)\sin nt\intl_0^t\overline{u^2(s)}\sin(2ns)dsdt-
\intl_0^\pi \frac{t}{\pi}f(t)\cos nt\intl_0^\pi
\overline{u(s)}\sin(2ns)dsdt+
$$\begin{equation} \label{8}
\left.+\intl_0^\pi\frac{tf(t)}{2\pi n}\cos nt\intl_0^\pi
\overline{u^2(s)}(\cos(2ns)-1)dsdt-
\frac{1}{2n}\intl_0^\pi f(t)\cos nt\intl_0^t
\overline{u^2(s)}\cos(2ns)dsdt+\right.
\end{equation}
$$
+\frac{1}{2n}\intl_0^\pi f(t)\cos nt\intl_0^t
\overline{u^2(s)}dsdt-\intl_0^\pi\frac{2tf(t)}{\pi}\cos
nt\intl_0^{\pi}\intl_0^s
\overline{u(s)}\overline{u(\tau)}\cos(2ns)\sin(2n\tau)d\tau dsdt+
$$
$$
+2\intl_0^\pi f(t)\cos nt\intl_0^t\intl_0^s
\overline{u(s)}\overline{u(\tau)}\cos(2ns)\sin(2n\tau)d\tau dsdt.
$$

Покажем, что для любой функции $f\in L_2[0,\pi]$ справедлива
оценка:
\begin{equation}
\label{Psi}\left\|\suml_{n=N}^{m}\Psi_n\sin nx\right\|_C\le
C_u\|f\|_{L_2}.
\end{equation}
Оценим первое слагаемое в полученной
сумме (для этого поменяем местами пределы интегрирования в двойном
интеграле):
\begin{multline}\label{9}
\left|\suml_{n=N}^m\frac{\sin nx}{2n}\intl_0^\pi f(t)\sin
nt\intl_0^t\overline{u^2(s)}\sin(2ns)dsdt\right|\le\\
\le\intl_0^\pi|u(s)|^2\left|\suml_{n=N}^m\frac{\sin(2ns)}{2n}\sin
nx\intl_s^\pi f(t)\sin ntdt\right|ds\le\\
\le\intl_0^\pi|u(s)|^2\suml_{n=N}^m\frac{1}{2n}\left|\intl_s^\pi
f(t)\sin ntdt\right|ds\le
 C_u \|f\|_{L_2},
\end{multline} так как последовательность $\left\{\intl_s^\pi f(t)\sin
ntdt\right\}_{n=N}^{\infty}=\left\{(H_sf,\sin
nt)\right\}_{n=N}^{\infty}$ принадлежит пространству $l_2$ и ее
норма в этом пространстве не превосходит
$\|H_sf\|_{L_2}\le\|f\|_{L_2}$ (здесь $H_s$ -- оператор срезки в
пространстве $L_2[0,\pi]$).

Абсолютно аналогичные рассуждения можно провести для четвертого слагаемого в правой части \eqref{8}.

Рассмотрим второе слагаемое.
$$ \left|\suml_{n=N}^m\intl_0^\pi \frac{t}{\pi}f(t)\cos
nt\intl_0^\pi \overline{u(s)}\sin(2ns)dsdt\sin nx\right|\le
C\left(\suml_{n=N}^m\left|\intl_0^\pi tf(t)\cos
ntdt\right|^2\right)^{1/2}\cdot
$$
\begin{equation}
\label{10}\cdot\left(\suml_{n=N}^m\left|\intl_0^\pi
\overline{u(s)}\sin(2ns)ds\right|^2\right)^{1/2}\le
C_u\|f\|_{L_2}.
\end{equation}

Третье слагаемое можно оценить похожим образом: $$
\left|\suml_{n=N}^m\intl_0^\pi\frac{tf(t)}{2\pi n}\cos
nt\intl_0^\pi \overline{u^2(s)}(\cos(2ns)-1)dsdt\sin nx\right|\le $$
\begin{equation}
\label{11} \le C_u
\left(\suml_{n=N}^m\frac{1}{n^2}\right)^{1/2}\left(\suml_{n=N}^m\left|\intl_0^\pi
tf(t)\cos ntdt\right|^2\right)^{1/2}\le C_u\|f\|_{L_2}.
\end{equation}

Пятое слагаемое в \eqref{8} преобразуем, проинтегрировав функцию
$\overline{u^2(s)}$: \begin{equation}
\label{11'}\left|\suml_{n=N}^m\frac{1}{2n}\intl_0^\pi f(t)\cos
nt\intl_0^t \overline{u^2(s)}dsdt\right|\le
\suml_{n=N}^m\frac{1}{2n}\left|\intl_0^\pi f(t)U(t)\cos nt
dt\right|\le C_u \|f\|_{L_2},\end{equation} поскольку функция
$U(t)=\intl_0^t\overline{u^2(s)}ds$ является абсолютно
непрерывной, следовательно, произведение $f(t)U(t)$ принадлежит
пространству $L_2[0,\pi]$.

Перейдем к тройным интегралам в правой части \eqref{8}. Рассмотрим
первый из них.
\begin{equation}
\label{12} \left|\suml_{n=N}^m\intl_0^\pi\frac{2tf(t)}{\pi}\cos
nt\intl_0^{\pi}\intl_0^s \overline{u(s)}\overline{u(\tau)}\cos(2ns)\sin(2n\tau)d\tau dsdt\sin nx\right|\le
\end{equation}
$$
 \le C\left(\suml_{n=N}^m\left|\intl_0^\pi tf(t)\cos
ntdt\right|^2\right)^{1/2}\left(\suml_{n=N}^m\left|\intl_0^\pi\intl_0^s
\overline{u(s)}\cos(2ns)\overline{u(\tau)}\sin(2n\tau)d\tau
ds\right|^2\right)^{1/2}\le
$$
$$
\le C\|f\|_{L_2}\left(\suml_{n=N}^m\left|\Big(u(s)u(\tau)\chi
(s,\tau),\cos(2ns)\sin(2n\tau)\Big)\right|^2\right)^{1/2}\le
C_u\|f\|_{L_2},
$$
поскольку функция $u(s)u(\tau)\chi(s,\tau)\in L_2[0,\pi]^2$

Наконец,
$$
\left|\suml_{n=N}^m\intl_0^\pi f(t)\cos nt\intl_0^t\intl_0^s
\overline{u(s)}\overline{u(\tau)}\cos(2ns)\sin(2n\tau)d\tau dsdt\sin nx\right|\le
$$
\begin{equation}
\label{13} \le\left|\intl_0^\pi u(s)\suml_{n=N}^m\cos(2ns)\sin
nx(\widetilde{H}_sf(t),\cos nt)
(H_su(\tau),\sin(2n\tau))ds\right|\le
\end{equation}
$$
 \le\intl_0^\pi|u(s)|\left|\suml_{n=1}^m(\widetilde{H}_sf(t),\cos nt)(H_su(t),\sin(2nt))\right|ds
\le \intl_0^\pi|u(s)|\|\widetilde{H}_sf\|_{L_2}\|H_su\|_{L_2}ds\le
C_u\|f\|_{L_2}.
$$
Здесь через $\widetilde{H}_s$ обозначен оператор срезки
$\widetilde{H}_sf(t)=\chi_{[s,\pi]}f(t)$.

Из \eqref{9}--\eqref{13} следует неравенство \eqref{Psi} (а
значит, и \eqref{0'}).

\medskip
{\bf  Шаг 3 (оценка третьего и четвертого слагаемых в основном
равенстве).}
\medskip

Оценка третьего слагаемого в \eqref{3} получается тривиально:
\begin{equation}
\label{14}\suml_{n=N}^m\left|\intl_0^\pi
f(t)\overline{\psi_{n,2}(t)}dt\right|\le C_u\|f\|_{L_2},
\end{equation}
так как последовательность $\{\|\psi_{n,2}(t)\|_C\}\in l_1$.

Перейдем к четвертому члену представления \eqref{3}. В силу асимптотических формул \eqref{efas}:
\begin{equation} \label{15}
\left\|\sum_{n=N}^{m}\sqrt{\frac{2}{\pi}}(f(t),\sin
nt)\varphi_n(x)\right\|_{C}\le
\sum_{n=N}^{m}(|f_n|\cdot\|\varphi_n(x)\|_{C})\le C_u\|f\|_{L_2},
\end{equation}
где $f_n=\sqrt{2/\pi}(f(x),\sin nx)$.

Наконец,
\begin{equation}
\label{16}
\left\|\sum_{n=N}^{m}(f(t),\psi_n(t))\varphi_n(x)\right\|_{C}\le
\|f\|_{L_2}\sum_{n=N}^{m}(\|\psi_n(t)\|_{L_2}\cdot\|\varphi_n(x)\|_{C})\le
C_u\|f\|_{L_2}.
\end{equation}
Из неравенств \eqref{0}, \eqref{0'} и \eqref{14}--\eqref{16} вытекает оценка \eqref{B}.

Теперь мы можем легко получить оценку для нормы оператора $B_m$.
Очевидно, что он представляется в виде суммы:
\(B_m=B_{N-1}+B_{m,N}\). Поскольку выбор числа \(N\) зависит
только от  \(u\), то ясно, что $\|B_{N-1}\|_{L_2\to C}\le C_u$
(так как этот оператор представляется в виде суммы конечного числа
слагаемых). Отсюда и из неравенства \eqref{B} немедленно следует,
что
\begin{equation} \label{B1}
\|B_{m}(u)\|_{L_2\to C}\le C_u.
\end{equation}

\medskip
{\bf Шаг 4 (доказательство равносходимости).}
\medskip

{\it Для любой функции \(f\in L_2[0,\pi]\) выполнено:
\begin{equation}\label{B0} \lim_{m\to\infty} \|B_mf\|_C=0. \end{equation}}

Рассмотрим действие оператора \(B_m\) на собственные и
присоединенные функции оператора \(L\):
$$
B_my_k(x)=\suml_{n=1}^{m}(y_k(x),w_n(x))y_n(x)-\frac{2}{\pi}\suml_{n=1}
^{m}(y_k(x),\sin nx)\sin nx.
$$
Первое слагаемое в правой части последнего соотношения равно 0 при
\(m<k\) и равно \(y_k(x)\) при \(m\ge k\). Второе слагаемое
представляет собой частичную сумму ряда Фурье функции \(y_k\). Так
как все функции \(y_k\in W_2^1 [0,\pi]\), то ряд Фурье функции
\(y_k\) сходится к ней равномерно на отрезке \([0,\pi]\), и мы
получаем, что $\liml_{m\to\infty}\|B_my_k\|_C=0. $

Осталось заметить, что, в силу полноты системы \(\{y_k(x)\}\) (см.
утверждение 1), из непрерывности оператора \(B_m\) следует
предельное соотношение \eqref{B0}.

Шаг 4 завершен.

\medskip

Перейдем к доказательству утверждения о скорости равносходимости.
Для любого $k\ge N_u$ обозначим
\(g_k(x)=\suml_{n=1}^{k}c_ny_n(x)\) (напомним, что
\(c_n=(f(x),w_n(x))\)). Очевидно, что для любой функции \(f\in
L_2[0,\pi]\) и для любого натурального \(m\) выполнено:
\begin{equation}
\label{17} \|B_m f\|_{C}\le\|B_m(f-g_k)\|_{C}+\|B_mg_k\|_{C}.
\end{equation}

\medskip
{\bf Шаг 5 (оценка нормы $B_m(f-g_k)$ в пространстве $C[0,\pi]$).}
\medskip

{\it Пусть \(g_k(x)=\suml_{n=1}^{k}c_ny_n(x)\), $k\ge N_u$. Тогда
\begin{equation} \label{18} \|B_m(f-g_k) \|_{C}\le
C_u\left(\suml_{n=k+1}^{\infty}|c_{n,0}|^2\right)^{1/2}+\|f\|_{L_2}\upsilon_u(k),
\end{equation}
где $c_{n,0}=\sqrt{2/\pi}(f(x),\sin nx)$, а $\upsilon_u(k)\to0$
при $k\to +\infty$.}

С учетом асимптотических формул \eqref{efas} получаем:
$$
\|f(x)-g_k(x)\|_{L_2}\le
\left\|\suml_{n=k+1}^{\infty}\frac{2}{\pi}(f(x),\sin nx)\sin
nx\right\|_{L_2}+
\left\|\suml_{n=k+1}^{\infty}\sqrt{\frac{2}{\pi}}(f(x),
\psi_n(x))\sin nx\right\|_{L_2}+$$
$$+\left\|\suml_{n=k+1}^{\infty}\sqrt{\frac{2}{\pi}}(f(x),\sin
nx)\varphi_n(x)\right\|_{L_2}+ \left\|\suml_{n=k+1}^{\infty}(f(x),
\psi_n(x))\varphi_n(x)\right\|_{L_2}\le$$
$$\le C_u
\left(\suml_{n=k+1}^{\infty}|c_{n,0}|^2\right)^{1/2}+\|f\|_{L_2}\upsilon_u(k),
$$
где $\upsilon_u(k)\to0$ при $k\to +\infty$.

В силу оценки \eqref{B1} на норму оператора $B_m$ из доказанного
неравенства немедленно вытекает неравенство \eqref{18}. Шаг 5
завершен.

Перейдем к оценке второго слагаемого в \eqref{17}. Пусть $m>k$.
Обозначим через $S_m$ оператор, действующий из пространства
$W_2^1[0,\pi]$ в пространство $C[0,\pi]$ по правилу:
$S_mh(x)=2/\pi\suml_{n=m+1}^{\infty}(h(t),\sin nt)\sin nx$.
Заметим, что
$$B_mg_k(x)=g_k(x)-\frac{2}{\pi}\suml_{n=1}^{m}(g_k(t),\sin nt)\sin
nx=S_mg_k(x)$$ (поскольку все собственные и присоединенные функции
оператора $L$ принадлежат пространству $W_2^1[0,\pi]$, то действие
оператора $S_m$ на них корректно определено). Тогда
\begin{equation}
\label{20} \|B_{m} g_k\|_{C}\le
\|S_m(g_k-g_{N})\|_{C}+\|S_mg_{N}\|_{C}\end{equation}

{\bf Шаг 6 (оценка нормы $S_m(g_k-g_{N})$ в пространстве
$C[0,\pi]$).}

{\it Пусть $k\ge N_u$ -- натуральное число. Тогда для любой
функции $f$ из пространства $L_2[0,\pi]$ и для любого $\eps>0$
справедлива оценка
\begin{equation}
\label{S}  \|S_m(g_k-g_{N})\|_{C}\le
C_{u,\eps}\|f\|_{L_2}km^{\eps-1/2}.
\end{equation}
}

Пусть $m>k\ge N_u$. Заметим, что левую часть неравенства \eqref{S} можно представить в виде
$$\|S_m(g_k(x)-g_{N}(x))\|_{C}=\left\|\suml_{n=N+1}^{k}c_nS_my_n(x)\right\|_{C}=
\left\|\suml_{n=N+1}^{k}c_nS_m\varphi_n(x)\right\|_{C},$$ так как
$y_n(x)=\sqrt{2/\pi}\sin nx+\varphi_n(x)$, где \(\varphi_n(x)\)
определены в \eqref{efas}. Далее,
$$
\left\|\suml_{n=N+1}^{k}c_nS_m\varphi_n(x)\right\|_{C}
\le\left(\suml_{n=N+1}^{k}|c_n|^2\right)^{1/2}
\left(\suml_{n=N+1}^{k}\|S_m\varphi_n(x)\|_C^2\right)^{1/2}\le
$$ $$\le
C_u\|f\|_{L_2}\left(\suml_{n=N+1}^{k}\|S_m\varphi_n(x)\|_{C}^2\right)^{1/2}.
$$
При выводе последнего неравенства мы учли тот факт, что
$$\suml_{n=N+1}^{k}|c_n|^2\le \suml_{n=N+1}^{\infty}|c_n|^2\le
2\left(\frac{2}{\pi}\suml_{n=N+1}^{\infty}|(f(x),\sin nx)|^2
+\suml_{n=N+1}^{\infty}|(f(x),\psi_n(x))|^2\right)\le
C_u\|f\|_{L_2}^2. $$

Оценим  $\|S_m\varphi_n(x)\|_C$, воспользовавшись теоремой
вложения Соболева: пространство
$W_2^{\varepsilon+1/2}[0,\pi]\hookrightarrow C[0,\pi]$ (\cite[п.
4.6.2]{Tr}). Получим, что для любого \(\varepsilon>0\),
$$
\|S_m\varphi_n(x)\|_{C}\le
C_{\varepsilon}\left(\suml_{j=m+1}^{\infty}(j^{2\varepsilon+1}+1)|(\varphi_n(x),\sin
jx)|^2\right)^{1/2}=
$$
\begin{equation}
\label{21}
=C_{\varepsilon}\left(\suml_{j=m+1}^{\infty}\frac{j^{2\varepsilon+1}+1}{j^2}|(\varphi'_n(x),\cos
jx)|^2\right)^{1/2}\le
C_{\varepsilon}m^{\varepsilon-1/2}\|\varphi_n(x)\|_{W_2^1}.
\end{equation}
Здесь мы воспользовались интегрированием по частям и учли тот
факт, что $\varphi_n(0)=\varphi_n(\pi)=0$.

Так как из асимптотических формул \eqref{efas} следует, что
$\varphi'_n(x)=n\eta_n(x)+u(x)y_n(x)$, то
\(\|\varphi_n(x)\|_{W_2^1}\le C_{u}n\eta_n\), где
\(\|\{\eta_n\}\|_{l_2}\le C_{u}\). Значит,
$$\|S_m\varphi_n(x)\|_{C}\le
C_{u,\varepsilon}m^{\varepsilon-1/2}n\eta_n.$$ Итак, первое
слагаемое в \eqref{20} можно оценить следующим образом: $$
\|S_m(g_k-g_{N})\|_{C}\le
C_{u,\varepsilon}\|f\|_{L_2}\left(\suml_{n=1}^{k}m^{2\varepsilon-1}n^{2}\eta_n^2\right)^{1/2}\le
C_{u,\varepsilon}\|f\|_{L_2}m^{\varepsilon-1/2}k.
$$

Шаг 6 завершен.

Перейдем ко второму слагаемому в \eqref{20}. Действие оператора
$S_m$ на функции $g_{N}$ оценивается точно так же, как в
\eqref{21}:
$$
\|S_mg_{N}(x)\|_{C}\le
C_{\varepsilon}m^{\varepsilon-1/2}\|g_{N}(x)\|_{W_2^1}.
$$
Так как \(\|g_{N}\|_{W_2^1}\le C_u\|f\|_{L_2}\), то
\begin{equation}
\label{S1}  \|S_mg_{N}\|_{C}\le
C_{u,\varepsilon}\|f\|_{L_2}m^{\varepsilon-1/2}.
\end{equation}
Из соотношений \eqref{20}, \eqref{S} и \eqref{S1} теперь сразу вытекает, что для любого натурального $k\ge N_u$
\begin{equation}
\label{22}\|B_mg_k\|_{C}=\|f\|_{L_2}\upsilon_u(m),
\end{equation}
где $\upsilon_u(m)\to 0$ при $m\to+\infty$. Оценка \eqref{1}
теперь следует из \eqref{17}, \eqref{18} и \eqref{22}. Теорема
полностью доказана.

\vskip 2cm

{\bf \S 3. Случай вещественного потенциала.}

\medskip

В этом разделе мы хотели бы анонсировать результат, являющийся в
случае вещественного потенциала несколько более общим. Дело в том,
что доказательство утверждения 2, приведенное в работе
\cite{SSh2}, опиралось на асимптотические формулы для собственных
значений оператора $L$, полученные в той же статье. Однако позже
теми же авторами были выведены асимптотики собственных значений с
оценками остатков, равномерными по шару \(u\in B_{R}=\{v\in
L_2[0,\pi]:\|v\|_{L_2}\le R\}\) (эти результаты еще не
опубликованы). С учетом полученных асимптотических формул,
рассуждениями, полностью аналогичными доказательству теоремы 3.13
работы \cite{SSh2}, можно получить формулы
\eqref{efas}--\eqref{efas1}, но уже с оценками норм всех
последовательностей, зависящими не от конкретного потенциала а
лишь от радиуса шара $R$.

В этом случае будет справедлива следующая

\begin{Theorem}
Пусть в условиях теоремы 1 функция $u$ является
вещественнозначной, причем $\|u\|_{L_2}\le R$ для некоторого
$R>0$. Тогда
\begin{equation}
\label{1R} \left\|\sum_{n=1}^{m}{c_n
y_n(x)}-\sum_{n=1}^{m}\sqrt{\frac{2}{\pi}}{c_{n,0}\sin
nx}\right\|_{C}\le
C_R\left(\suml_{n\ge\sqrt{m}}^{\infty}|c_{n,0}|^2\right)^{1/2}+\|f\|_{L_2}\upsilon_R(m),
\end{equation}
где $\upsilon_R(m)\to 0$ при $m\to+\infty$.
\end{Theorem}

Автор благодарит проф. А.А.Шкаликова и доц. А.М.Савчука за
полезные замечания.

\medskip

\end{document}